%%%%%%%%%%%%%%% written by AMS-TEX %%%%
%%%%%%%%%%%%%%%%%%

\input amstex
\documentstyle{amsppt}
\NoRunningHeads
%\NoPageNumbers
\magnification = \magstep 1
%\NoBlackBoxes
\pagewidth{32pc}
\pageheight{48pc}

\topmatter
\title\  The equality
$I^2=QI$ in Buchsbaum rings with
multiplicity two
\endtitle

\abstract Let $A$ be a Buchsbaum local
ring with the maximal ideal $\frak{m}$
and let $\roman{e}(A)$ denote the
multiplicity of $A$. Let $Q$ be a
parameter ideal in $A$ and put $I = Q
: \frak{m}$. Then the equality $I^2=QI$
holds true, if
$\roman{e}(A) = 2$ and
$\roman{depth}~A > 0$.  The assertion
is no longer true, unless
$\roman{e}(A) = 2$. Counterexamples
are given.
\endabstract

\leftheadtext{Shiro Goto and Hideto
Sakurai}
\rightheadtext{The equality
$I^2=QI$ in Buchsbaum local rings with
multiplicity two}

\author Shiro Goto and Hideto
Sakurai
\endauthor

\address Department of Mathematics, School of Science and
Technology, Meiji University, 214-8571 JAPAN 
\endaddress  

\email goto\@math.meiji.ac.jp \endemail

\address Department of Mathematics, School of Science and
Technology, Meiji University, 214-8571 JAPAN 
\endaddress  

\email ee78052\@math.meiji.ac.jp
\endemail

 

\thanks The first author is supported
by the Grant-in-Aid for Scientific
Researches in Japan (C(2), No. 13640044).
\endthanks
\endtopmatter

\def\depth{\operatorname{depth}}

\def\Ass{\operatorname{Ass}}

\def\Min{\operatorname{Min}}

\def\bar{\overline}

\def\m{\frak {m}}

\def\K{\roman{K}}

\def\e{\roman{e}}

               %reduction number
              %type
\def\l{\ell}
\def\n{{\frak n}}

\def\p{{\frak p}}
\def\p{{\frak p}}
\def\q{{\frak q}}
\def\H{{\roman H}}

\def\z{{\Bbb Z}}

\def\z{{\Bbb Z}}

\def\depth{\operatorname{depth}}

\def\Ass{\operatorname{Ass}}

\def\Min{\operatorname{Min}}

\def\bar{\overline}

\def\l{\lambda}

\def\m{\frak m}
\def\n{\frak n}
\def\q{\frak q}
\def\p{{\frak p}}

               %reduction number
\def\H{\roman H}

\def\e{\roman e}

\font\b=cmr10 scaled \magstep4
%\font\c=cmr10 scaled \magstep3

\def\bigzerou{\smash{\lower1.7ex\hbox{\b 0}}}

\document
\baselineskip=15pt

%%%%%%%%%%%%%%%%%%%%%%%%%           1        %%%%%%%%%%%%%%%%%%%%%%%%%%

\head 1. Introduction.
\endhead 

Let $A$ be a Noetherian local
ring with the maximal ideal $\frak{m}$
and $d =
\dim A$. Let $Q$ be a parameter ideal
in $A$ and let $I = Q : \frak{m}$.
In this paper we are interested in the
problem of when the equality
$I^2 = QI$ holds true. This problem was
completely solved by A. Corso, C.
Huneke, C. Polini, and W. Vasconcelos
\cite{CHV, CP, CPV} in the case where
$A$ is a Cohen-Macaulay ring. 
When
$A$ is a Buchsbaum ring, partial
answers only recently appeared in
the authors' paper \cite{GSa},
supplying 
\cite{Y1, Y2} and \cite{GN} with ample
examples of ideals $I$, for which the
Rees algebras
$\roman{R}(I)=
\bigoplus_{n
\geq 0}I^n$, the associated graded
rings $\roman{G}(I) =
\roman{R}(I)/I\roman{R}(I)$, and the
fiber cones
$\roman{F}(I)=\roman{R}(I)/\frak{m}\roman{R}(I)$
are all Buchsbaum rings with certain
specific graded local cohomology
modules.

This research is
a succession of \cite{GSa} and the 
present  purpose is to prove the
following, in which $\roman{e}(A) =
\roman{e}_{\frak{m}}^0(A)$ denotes the
multiplicity of $A$ with respect to the maximal ideal 
$\frak{m}$.

\proclaim{Theorem (1.1)}
Let $A$ be a Buchsbaum ring.
Assume that
$\roman{e}(A) = 2$ and $\depth A > 0$.
Then the equality $I^2 = QI$ holds true
for all parameter ideals $Q$ in $A$,
where $I = Q :\frak{m}$.
\endproclaim

\noindent
The readers may consult \cite{G2} for the structure
of Buchsbaum local rings $A$ with $\e(A) = 2$. There is
given in \cite{G2, Sections 3,4} a complete list of
equi-characteristic non-Cohen-Macaulay Buchsbaum complete 
local rings $A$ with $\e(A) = 2$, $\depth A > 0$, and
infinite residue class fields.

In their remarkable papers \cite{CP,
CPV} A. Corso, C. Polini, and
W. Vasconcelos proved that the
equality $I^2 = QI$ holds true for
every parameter ideal $Q$ in
$A$, if $A$ is a Cohen-Macaulay local
ring with
$\roman{e}(A) \geq 2$. 
This is no more true in the Buchsbaum
case. As is shown
in \cite{GSa, Theorem
(4.8)}, for every
integer $e \geq 3$, there exists a
Buchsbaum local ring $A$ with $\dim A
= 1$ and $\roman{e}(A) = e$ which
contains a parameter ideal $Q$ such
that $I^e = QI^{e-1}$ but $I^{e-1} \ne
QI^{e-2}$. Accordingly, without
additional assumptions on Buchsbaum local
rings $A$, no hope is left for the
equality
$I^2 = QI$, at least in the case where
$\dim A= 1$ and
$\roman{e}(A) \geq 3$. Our theorem 
(1.1) settles the
case where
$\roman{e}(A) = 2$ and $\depth A > 0$,
providing a drastic break-through
against the counter-examples of
\cite{GSa}.

Before entering details, let us
briefly note how this paper is
organized. We shall prove Theorem
(1.1) in Section 3. For the purpose
we need some preliminaries, that we
will summarize in Section 2. The
counterexamples given by
\cite{GSa} are all of dimension 1,
and according to Theorem (1.1), it
might be natural to suspect that
the equality $I^2 = QI$ holds true in
higher dimensional cases of higher
depth. The answer is, nevertheless,
still negative. We shall construct
examples, showing that for
given integers $1
\leq d < m$, there exists a Buchsbaum
local ring
$A$ with $\dim A = d$, $\depth A =
d-1$, and
$\roman{e}(A) = 2m$, which contains a
parameter ideal $Q$ such that $I^3 =
QI^2$ but $I^2 \ne QI$ (Theorem (4.5) and
Proposition (4.7)).

\head 
2. Preliminaries.
\endhead

The purpose of this section is to
summarize some preliminary steps,
which we need to prove Theorem (1.1).
The result might have its own
significance. In such a case we shall
include a closer proof.

Here let us
fix our standard notation. Otherwise
specified, let $A$ be a Noetherian local
ring with the maximal ideal $\frak{m}$
and $d = \dim A$. For an ideal
$\frak{a}$ let
$\frak{a}^{\sharp}$ be the integral
closure of $\frak{a}$. Let $\ell_A(*)$
and $\mu_A(*)$ respectively denote the
length and the number of generators.
When $A$ is a Cohen-Macaulay local ring, we
denote by $\roman{r}(A)$ the
Cohen-Macaulay type of $A$, that is
$$\roman{r}(A) =
\ell_A(\roman{Ext}_A^d(A/\frak{m},A)).$$
Let $\roman{H}_{\frak{m}}^i(*)$ $(i
\in \z)$ be the local cohomology
functors of $A$ with respect to
$\frak{m}$. We denote by $\roman{e}(A)
=
\roman{e}_{\frak{m}}^0(A)$  the
multiplicity of
$A$.

Let us begin with the following.

\proclaim{Lemma (2.1)}
Assume that $\roman{e}(A) \geq 2$. Then
$\mu_A(I) =
\ell_A(I/Q) + d$
for every parameter ideal $Q$ in
$A$, where $I = Q : \frak{m}$. 
\endproclaim

\demo{Proof}
By \cite{GSa, Proposition (2.3)} $Q$ is
a minimal reduction of $I$, whence
$\frak{m}I = \frak{m}Q$. Thus
$$\mu_A(I) = \ell_A(I/\frak{m}I) =
\ell_A(I/Q) + \ell_A(Q/\frak{m}Q) =
\ell_A(I/Q) + d$$
as is claimed.
\qed
\enddemo

\proclaim{Proposition (2.2)}
Suppose that $A$ is a Cohen-Macaulay
local ring with $d = \dim A\geq 1$ and
let
$Q$ be a parameter ideal in $A$. Then 
$$
\ell_A((0) :_{\frak{m}/Q\frak{m}}
\frak{m}) =
\cases \roman{r}(A) + d & \text{if} 
\quad  Q
\ne Q^{\sharp},\\
d &  \text{if} \quad  Q =
Q^{\sharp}.
\endcases$$
\endproclaim

\demo{Proof}
Since $\dim A \geq 1$, we have
$Q\frak{m} : \frak{m} \subseteq
\frak{m}$, whence
$$\ell_A((0) :_{\frak{m}/Q\frak{m}}
\frak{m}) = \ell_A((Q\frak{m} :
\frak{m})/Q\frak{m}).$$ Let $I = Q :
\frak{m}$. Firstly, assume that $Q \ne
Q^{\sharp}$. Then $Q$ is a minimal
reduction of $I$, because $I^2 = QI$ (cf.
e.g., \cite{GH, Proposition (3.4)}).
Hence $\frak{m}I = \frak{m}Q$, so that
we have $I
\subseteq Q\frak{m} : \frak{m}$; thus $I
= Q\frak{m} : \frak{m}$. Consequently,
$$\ell_A((0) :_{\frak{m}/Q\frak{m}}
\frak{m}) = \ell_A(I/\frak{m}Q) =
\ell_A(I/Q) + \ell_A(Q/\frak{m}Q) =
\roman{r}(A) + d.$$ Suppose that $Q =
Q^{\sharp}$. Then $A$ is a regular
local ring which contains a regular
system $a_1, a_2, \cdots, a_d$ of
parameters such that $Q = (a_1,
\cdots, a_{d-1}, a_d^q)$ for some $q \geq
1$ (\cite{G3, Theorem (3.1)}). Hence
$Q\frak{m} : \frak{m} = Q$, because $\m = (a_1, a_2,
\cdots, a_d)$ and $Q\m : \m \subseteq (a_1,
\cdots, a_{d - 1}, a_d^{q+1}) :a_d =Q$. Thus
$$\ell_A((0) :_{\frak{m}/Q\frak{m}}
\frak{m}) = \ell_A(Q/\frak{m}Q) = d$$
as is claimed.
\qed
\enddemo

Let $A$ be a Buchsbaum local ring with
the Buchsbaum invariant $\roman{I}(A)$.
Then all the local cohomology modules
$\roman{H}_{\frak{m}}^i(A)$ $(i \ne d)$
are killed by the maximal ideal
$\frak{m}$ and one has the equality
$$\roman{I}(A) =
\sum_{i=0}^{d-1}{d-1\choose i}h^i(A),$$
where
$h^i(A)=\ell_A(\roman{H}_{\frak{m}}^i
(A))$ (\cite{SV, Chap.~I, Proposition
2.6}). Let
$$\roman{r}(A) =
\underset{Q}\to{\sup}~
\ell_A((Q:\frak{m})/Q)$$
where
$Q$ runs over parameter ideals in $A$ 
and call it the Cohen-Macaulay type
of $A$. We then have
$$\roman{r}(A) =
\sum_{i=0}^{d-1}{d\choose i}h^i(A) +
\mu_{\hat{A}}(\roman{K}_{\hat{A}})\tag
2.3$$ (\cite{GSu, Theorem (2.5)}),
where
$\roman{K}_{\hat{A}}$ denotes the
canonical module of the $\frak{m}$-adic
completion $\widehat{A}$ of $A$. Consequently $\roman{r}(A) <
\infty$.

\proclaim{Theorem (2.4)}
Let $B$ be a Gorenstein local ring with
$d = \dim B \geq 2$. Let $A$ be a
subring of $B$. Assume that $B$ is
a module-finite extension of $A$ and 
$\ell_A(B/A) = 1$. Then 
\roster
\item $A$ is a
Buchsbaum local ring and
$\dim A = \roman{I}(A) = d$. 
\item The equality $I^2 = QI$ holds
true for all parameter ideals $Q$ in
$A$, where $I = Q : \frak{m}$.
\endroster
\endproclaim

\demo{Proof}
Since $B$ is a module-finite extension
of $A$, by Eakin-Nagata's theorem our
ring $A$ is a Noetherian local ring with
$d = \dim~A$. Let $\frak{m}$ and
$\frak{n}$ be the maximal ideals in $A$
and $B$, respectively. We look at the
exact sequence
$$0 \to A \overset{\iota}\to{\to} B \to
A/\frak{m} \to 0 \tag 2.5$$
of $A$-modules, where $\iota$ denotes
the inclusion map. Then applying
functors $\roman{H}_{\frak{m}}^i(*)$ to
(2.5), we get
that $$\roman{H}_{\frak{m}}^i(A) =
(0)~(i
\ne 1, d) \ \  \text{and} \ \
\roman{H}_{\frak{m}}^1(A) \cong
A/\frak{m},$$ because $\depth_AB = d$.
Hence $A$ is a Buchsbaum ring with
$\roman{I}(A) = d$ (cf. \cite{SV,
Chap.~I, Proposition 2.12}). Notice
that
$\frak{m}B =
\frak{m}$, because $\frak{m}{\cdot}(B/A)
= (0)$. Hence $\frak{m}$ is an ideal
in $B$. On the other hand, we
naturally have by (2.5)  the exact
sequence $$0
\to A/\frak{m} \to B/\frak{m}B \to
A/\frak{m}
\to 0,$$ whence $\mu_A(B) = 2$. Let us
write
$B = A + At$ with $t \in  B$. Then $t
\not\in A$. We have $\roman{r}(A) = 
d + 2$ by (2.3), since $\roman{K}_A=
B$. Notice that
$\roman{e}(A) \geq 2$, because
$A$ is not a regular local ring.

Now let $Q = (a_1, a_2, \cdots, a_d)$ be
a parameter ideal in $A$ and put $I = Q
: \frak{m}$. Then $I$ is an ideal in
$B$, since so is $\frak{m}$. Thus $QB
\subseteq I$. We need the following.

\proclaim{Claim (2.6)}
$\ell_A(QB/Q) = d$.
\endproclaim

\demo{Proof of Claim (2.6)}
Since $B = A + At$,  we have $QB = Q +
\sum_{i = 1}^dA{\cdot}a_it$. Let
$\bar{a_it}$ denote the reduction of
$a_it$ mod $Q$. Then $QB/Q =
\sum_{i=1}^dk{\cdot}\bar{a_it}$\ \ ($k
= A/\frak{m}$). Let $\alpha_i \in A$
$(1
\leq i \leq d)$ and assume that
$\sum_{i=1}^d \alpha_i(a_it) \in Q$.
We write
$\sum_{i=1}^d \alpha_i(a_it) =
\sum_{i=1}^d \beta_ia_i$ with $\beta_i
\in A$. Then
$\sum_{i=1}^da_i(\alpha_it -
\beta_i)=0$. Because $a_1, a_2,
\cdots, a_d$ forms a
$B$-regular sequence, $\alpha_it -
\beta_i \in (a_j \mid j \ne i)B
\subseteq A$, so that $\alpha_it \in
A$. Hence
$\alpha_i \in \frak{m}$, because $t
\not\in A$. Thus the classes
$\{\bar{a_it}\}_{1
\leq i \leq d}$ form a $k$-basis of
$QB/Q$. Hence $\ell_A(QB/Q) = d$. 
\qed
\enddemo

If $\ell_A(I/Q) = \roman{r}(A)$, then
$I^2 = QI$ by \cite{GSa, Theorem
(3.9)}. Therefore to prove $I^2 =
QI$, we may assume  that  $\ell_A(I/Q)
\leq d + 1$. Hence, either  $\ell_A(I/Q)
= d$, or
$\ell_A(I/Q) = d + 1$  (cf.
Claim (2.6)). If
$\ell_A(I/Q) = d$, then $I = QB$, so
that
$I^2 = QB{\cdot}IB = QI$. Assume that
$\ell_A(I/Q) = d+1$. Then $\ell_A(I/QB) =
1$. Therefore
$\ell_B(I/QB) = 1$ and 
$\frak{n} I
\subseteq QB$. Hence $I = QB :
\frak{n}$, because
$QB \subsetneq I \subseteq QB :
\frak{n}$ and $B/QB$ is an Artinian
Gorenstein local ring. Accordingly, $I^2
= QB{\cdot}IB = QI$, if $QB \ne
(QB)^{\sharp}$ in $B$ (cf., e.g., \cite{GH,
Proposition (3.4)}). Suppose that $QB =
(QB)^{\sharp}$ in $B$. Then, since $\roman{e}(A) \geq 2$, by \cite{GSa,
Proposition
(2.3)} we have $I \subseteq Q^{\sharp}$. Hence $I \subseteq 
(QB)^{\sharp} = QB$ so that $I = QB$, which is
impossible, because $\ell_A(I/QB) =
1$. Thus
$QB
\ne (QB)^{\sharp}$ in $B$ and $I^2 =
QI$, which completes the proof of
Theorem (2.4). 
\qed
\enddemo

The proof of the following result
(2.7) is essentially the same as that
of Theorem (2.4). Let us indicate a
sketch only.

\proclaim{Proposition (2.7)}
Let $B$ be a Gorenstein local ring with
the maximal ideal $\frak{n}$ and $d =
\dim B
\geq 2$. Let $A$ be a subring of $B$
such that $B$ is a finitely generated
$A$-module. Assume that $A \subsetneq B$
and $\frak{n} \subseteq A$. Then 
\roster
\item $A$ is
a Buchsbaum local ring with $\frak{n}$ the maximal
ideal and
$\roman{I}(A) = d{\cdot}\ell_A(B/A)$.
\item The equality $I^2 = QI$ holds
true for all parameter ideals $Q$ in
$A$, where
$I = Q : \frak{m}$. 
\endroster
\endproclaim

\demo{Proof}
Similarly as in the proof of Theorem
(2.4), $A$ is a Buchsbaum local ring
with $\n$ the maximal ideal,
$\H_{\m}^i(A) = (0)$ ($i
\ne 1, d$), and $\H_{\m}^1(A) \cong
B/A$.  Let $Q$ be a parameter
ideal in $A$ and put $I = Q : \m$. 
Then $QB \subseteq I
\subseteq  QB :
\n$, since $I$ is an ideal of $B$ 
and $\m = \n$ in our case.
Therefore, either $QB = I$, or $I =
QB:\n$, since $B/QB$ is an Artinian
Gorenstein local ring. We
certainly have
$I^2=QI$ if $QB = I$. Assume that $I
\ne QB$.  Then $I = QB:\n$ and
$I^2 = QB{\cdot}IB = QI$,  because $QB
\ne (QB)^{\#}$ in $B$ for the same
reason as in the proof of Theorem
(2.4).
\qed
\enddemo

Before closing this section let us note one example
satisfying the hypothesis of
Proposition (2.7).

%\proclaim{Example (2.8)} 
\example{Example (2.8)}
Let $K/k$ be a finite extension of fields and assume that 
$\delta = [K:k] \geq 2$.
Let $n = \delta - 1$ and choose a
$k$-basis 
$\{ \theta_0 = 1,
\theta_1,\cdots, \theta_n \}$ of $K$.
Let $d \geq 2$ be an integer and let 
$B = K[[ X_1, X_2, \cdots, X_d]]$ be
the formal power series ring over
$K$. Let
$$A = k[[ \theta_i X_j \mid 0 \leq i \leq n, 1 \leq j \leq d ]].$$
Then $B$ is a module-finite extension of $A$ such that the
maximal ideal $\n$ of $B$ coincides with that of $A$,
$\ell_A(B/A) = n$, and
$\e(A) = \delta $. Hence by 
Proposition $(2.7)$ $A$ is a Buchsbaum
local  ring, in which the equality
$I^2 = QI$ holds true for all
parameter ideals $Q$, where $I =
Q:\m$.
%\endproclaim
\endexample

\demo{Proof}
Let $\frak{m}$ be the maximal
ideal in $A$, that is $\frak{m} =
(\theta_iX_j
\mid 0 \leq i \leq n, 1 \leq j \leq
d)A.$ Then $\frak{m}B = \frak{n}$
whence
$B =
\sum_{i=0}^nA
\theta_i$.  Let
$b
\in B$ and write
$b =
\sum_{i=0}^n a_i
\theta_i$ with $a_i \in A$.
Then, since $bX_j = \sum_{i=0}^n a_i (\theta_i X_j) \in \m$
for all $1 \leq j \leq d$, we get 
$\n \subseteq \m$. Hence $\m = \n$.
We have $$\mu_A(B) = \ell_A(B/\m B) =
\ell_A(B/\n) = [K:k] = \delta \geq
2.$$ Because $\ell_A(B/A) =
\ell_A(B/\n) -
\ell_A(A/\m)$, we get $\ell_A(B/A) =
\delta -1 = n$. Let $\q
=(X_1,X_2,\cdots, X_d)A$. Then since
$\m^2 = \q \m$, the ideal $\q$ is a
minimal reduction of $\m$, so that we
have
$$\e(A) = \e_{\q}^0(A) = \e_{\q}^0(B) =
\ell_A(B/\q B) = \ell_A(B/\n)=
\delta,$$
as is claimed.
\qed
\enddemo

\head 3. Proof of Theorem (1.1) \endhead

Let $A$ be a Noetherian local ring with $d=\dim A \geq 2$.
Suppose that $A$ is a reduced ring with $\# \Ass A = 2$, say
$\Ass A = \{ \p_1, \p_2 \}$. We furthermore assume that
$A/\p_i$ is a regular local ring with $\dim A/{\p_i} = d$
($i=1,2$) and
$\p_1 + \p_2 = \m$. Hence $\m = \p_1 \oplus \p_2$, because
$\p_1 \cap \p_2 = (0)$. With this notation and assumption we
have the following.

\proclaim{Theorem (3.1)}
$I^2 = QI$ for all parameter ideals
$Q$ in $A$, where $I = Q:\m$.
\endproclaim 

\demo{Proof}
We look at the exact sequences
$$0 \to A \overset \iota\to\rightarrow A/\p_1 \oplus A/\p_2 \to A/\m 
\to 0 \ \ \ \text{and} \tag 3.2$$
$$0 \to \p_i \to A \to A/\p_i \to 0 \tag 3.3 $$
of $A$-modules ($i = 1,2$), where $\iota(a) = (a \mod \p_1,a \mod \p_2)$
for all $a \in A$. Then, applying 
functors $\H_{\m}^i(*)$ to (3.2),
we get $\H_{\m}^i(A)=(0)$ ($i \ne
1,d$) and $\H_{\m}^1(A) \cong A/\m$.
Hence $A$ is a Buchsbaum local ring
with $\roman{I}(A) = d$ and
$\roman{e}(A) = 2$. We have 
$\roman{r}(A) = d+2$  by
(2.3), because $\K_A
= A/\p_1
\oplus A/\p_2$. Let $Q =
(a_1,a_2,\cdots,a_d)$ be a parameter
ideal in $A$ and put $I = Q:\m$. Then
$\ell_A(I/Q) \leq r(A) = d+2$. We may assume
$\ell_A(I/Q) \leq d+1$ (cf. \cite{GSa, Theorem (3.9)}).
Let $A_i = A/\p_i$ and $\m_i = \m/\p_i$ ($i=1,2$).
We write $a_j = \ell_j + m_j$ ($1 \leq j \leq d$) with
$\ell_j
\in
\p_1$ and
$m_j \in \p_2$.

Firstly we consider the case $QA_2 \ne (QA_2)^{\#}$ in
$A_2$. Let $\varepsilon : A \to A_2$ be the canonical
epimorphism. Then $\varepsilon(\p_1) = \m_2$ and $\p_1 \cong
\m_2$ via $\varepsilon$, because $\m =
\p_1 \oplus \p_2$. Hence 
$$\ell_A((0):_{\p_1/Q \p_1} \m) = \ell_{A}((0):_{\m_2/Q
\m_2}
\m) = \ell_{A_2}((0):_{\m_2/Q \m_2} \m_2) = d+1 \tag 3.4 $$
by Proposition (2.2). Let $\bar{*}$ denote the reduction mod $\p_2$.
Then since $A_2/QA_2$ is an
Artinian Gorenstein local ring
and since $Q\m_2 :
\m_2 = QA_2 :
\m_2 \subseteq \m_2$  (cf. Proof of
Proposition (2.2)), the the ideal $Q
\m_2:\m_2$ of $A_2$ is generated by
$\{ \bar{\ell_j} \}_{1 \leq j
\leq d}$ together with one 
 more element, say $\bar{\eta}$\
\ ($\eta
\in
\p_1)$. Hence the $A/\m$-space $(0)
:_{\p_1/Q \p_1} \m$ is spanned by 
$\{\ell_j \mod Q\p_1 \}_{1 \leq j 
\leq d}$ and $\eta \mod
Q\p_1$, because $\p_1 \cong \m_2$ via $\varepsilon$.
Now look at the exact sequence
$$0 \to \p_i/Q\p_i \to A/Q \overset
\varphi_i\to\rightarrow A_i/QA_i \to 0
\ \ \ \ \text{($i=1,2$)}
\tag 3.5$$  
of $A$-modules induced from (3.3) (notice that
$a_1,a_2,\cdots,a_d$ is
an $A_i$-regular sequence). Then,
considering the socles of the terms
in  (3.5) with $i=1$, we get $$I= Q +
(\ell_1,\ell_2,
\cdots, \ell_d)+(\eta),$$ 
because
$\ell_A((0) :_{\p_1/Q\p_1} \m ) =
d+1$ by (3.4) and
$\ell_A(I/Q) \leq d+1$ by our standard assumption.
Hence $I^2 = QI + (\eta^2)$, because
$$\ell_i \ell_j = (\ell_i +
m_i)\ell_j = a_i \ell_j \ \ \
\text{and} \ \ \
\ell_i
\eta = (\ell_i + m_i) \eta = a_i
\eta$$ for all $1 \leq i,j \leq d$
\ (recall that $\p_1 \cap \p_2 =
(0)$). On the other hand, since $QA_2
\ne (QA_2)^{\#}$, we get
$(QA_2:\m_2)^2 = QA_2 \cdot
(QA_2:\m_2)$  (cf., e.g., \cite{GH,
Proposition (3.4)}).  Hence
$\bar{\eta}^2 \in
(\bar{\ell_1},\bar{\ell_2},\cdots,\bar{\ell_d}){\cdot}[(\bar{\ell_1},\bar{\ell_2},\cdots,\bar{\ell_d})+
(\bar{\eta})]$ so that 
$$\eta^2 \in
(\ell_1,\ell_2,\cdots,\ell_d){\cdot}[(\ell_1,\ell_2,\cdots,\ell_d)+(\eta)]
+ \p_2.
\tag 3.6$$
Because $\p_1 \cap \p_2 = (0)$ and $(\ell_1,\ell_2,\cdots,\ell_d)+(\eta)
\subseteq
\p_1$,
by (3.6) we readily get that 
$$\eta^2 \in
(\ell_1,\ell_2,\cdots,\ell_d){\cdot}[(\ell_1,\ell_2,\cdots,\ell_d)+(\eta)]
\subseteq QI.$$ Hence $I^2 = QI$,
since
$I^2= QI +(\eta^2)$. 
We get, by the symmetry, that $I^2 =
QI$ also in the
case where
$QA_1 \ne (QA_1)^{\#}$.

We now consider the case $QA_i = 
(QA_i)^{\#}$ for $i = 1,2$.
Then thanks to the exact sequence (3.5) with $i=1$, we have 
$\ell_A(I/Q) \geq d$, because
$$\ell_A((0):_{\p_1/Q\p_1} \m) = 
\ell_{A_2}((0):_{\m_2/Q\m_2} \m_2) =
d$$ by Proposition (2.2). We actually
have the following .

\proclaim{Claim (3.7)}
$\ell_A(I/Q)=d$.
\endproclaim

\demo{Proof of Claim (3.7)}
Suppose that $\ell_A(I/Q) \ne d$. Then $\ell_A(I/Q) = d+1$.
Choose a regular system $\bar{c_1},\bar{c_2},\cdots,\bar{c_d}$
($c_i \in \p_1$) of parameter for $A_2$ so that 
$QA_2 =(\bar{c_1},\cdots,\bar{c_{d-1}},\bar{c_d}^q)$ for
some 
$q > 0$ (cf. \cite{G3, Theorem
(3.1)}). Hence $$QA_2 : \m_2 =
(\bar{c_1},\cdots,\bar{c_{d-1}},\bar{c_d}^{q-1})
= QA_2 +  (\bar{c_d}^{q-1}).$$ We
have 
$$ (\ell_1,\ell_2,\cdots,\ell_d) =
(c_1,\cdots,c_{d-1}, c_d^q),
\tag 3.8 $$
because $QA_2 = (\bar{\ell_1},\bar{\ell_2},\cdots,\bar{\ell_d}) = 
(\bar{c_1},\cdots,\bar{c_{d-1}},\bar{c_d}^{q})$ 
and $\p_1 \cong \m_2$ via $\varepsilon$.
Notice that $$\ell_A((0):_{\p_2/Q\p_2} \m) = 
\ell_{A_1}((0):_{\m_1/Q\m_1} \m_1) = d$$ and 
$(0):_{\m_1/Q\m_1} \m_1 = QA_1/Q\m_1$
(cf. Proposition (2.2) and its proof).
Hence the $A/\m$-space $(0):_{\p_2/Q\p_2} \m$ is spanned by 
$\{ m_j \mod Q\p_2 \}_{1 \leq j \leq d}$.
We now look at the exact sequence
(3.5) with $i=2$. Then since
$\ell_A((0):_{\p_2/Q\p_2} \m) = d $
and 
$\ell_A(I/Q) = d+1$,
the canonical epimorphism $\varphi_2
: A/Q
\to A_2/QA_2$ cannot be zero on the
socles,  so that we have $IA_2 =
QA_2:\m_2.$ Hence $IA_2 = QA_2 +
(\bar{c_d}^{q-1})$. Choose $\eta \in
I$ so that $\bar{\eta} =
\bar{c_d}^{q-1}$ and write $\eta =
c_d^{q-1} + \delta + \rho$ with 
$\delta \in Q$ and $\rho \in \p_2$.
Then thanks to the exact sequence 
(3.5) with $i=2$, we get
$$
\align
I &= Q + (m_1,m_2,\cdots,m_d) +
(\eta) \\ &=Q + (m_1,m_2,\cdots,m_d)
+  (c_d^{q-1} + \delta + \rho),
\endalign$$
because
$(0):_{\p_2/Q\p_2} \m$ is spanned by $\{ m_j \mod Q\p_2 \}_{1 \leq j
\leq
d}$.
Consequently we have
$$\matrix
\format \l & \l & \l \\
I & =  Q + (m_1,m_2,\cdots,m_d) +
(c_d^{q-1} + \rho) & 
                             \ \ \ \text{(since $\delta \in
Q$)}\\ 
   & =  (\ell_1,\ell_2,\cdots,\ell_d) + (m_1,m_2,\cdots,m_d) +
(c_d^{q-1} + \rho)  & \\
  & =  (c_1,\cdots,c_{d-1},c_d^q) + (m_1,m_2,\cdots,m_d) + (c_d^{q-1} +
\rho)
&
                             \ \ \  \text{(by (3.8))}  \\
 & =  (c_1,\cdots,c_{d-1}) + (m_1,m_2,\cdots,m_d) + (c_d^{q-1} + \rho) &
\ \  \ \text{(since $c_d^{q} = c_d(c_d^{q-1}+ \rho)$)},
\endmatrix$$
which is impossible, because
$\mu_A(I) = \ell_A(I/Q) = 2d+1$ by
Lemma (2.1). Thus 
$\ell_A(I/Q) = d$. 
\qed
\enddemo

Therefore, in the exact sequence (3.5)
with
$i=2$, the socle $I/Q$ of $A/Q$
coincides with the image of $(0):_{\p
_2/Q\p _2} \m$, because
$\ell _A((0):_{\p _2/Q\p _2} \m)=d$;
that is $$I=Q+( m_1,m_2,\cdots,m_d).$$
Thus $I^2=QI$, because
$m_im_j=(\ell_i+m_i)m_j \in QI$ 
for all $1 \leq i,j \leq d$. This
completes the proof of Theorem (3.1).
\qed
\enddemo

We are now ready to prove Theorem(1.1).

\demo{Proof of Theorem (1.1)}
Passing to the ring $A[X]_{\m A[X]}$ 
where $X$ is an indeterminate over
$A$ and then passing to the
completion, we may assume that $A$ is
a complete local ring with the
infinite residue class field. Thanks
to [CP, Theorem 2.2], we may assume
that $A$ is not a Cohen-Macaulay ring.
Hence
$d \geq 2$ and so by \cite{G2, Theorem
1.1} 
$$\H_{\m} ^i(A)=(0)
\ ( i \ne 1, d ) \ \ \text{and} \
\ \H_{\m}^1(A)
\cong A/\m.$$ Let $B$ be the
Cohen-Macaulayfication of $A$, that
is the intermediate ring $A \subseteq
B \subseteq Q(A)$ where $Q(A)$
denotes the total quotient ring of
$A$, such that $B$ is a module-finite
extension of $A$,
$\depth_A B=d$, and $\m B = \m$ (cf.
\cite{G1, Theorem (1.1)}).  Then
$\H_{\m}^1(A) \cong B/A$ and hence
$\ell_A(B/A) = 1$. Now there are two
cases. One is the case where
$B$ is a local ring. The other one is
the case where
$B$ is not a local ring.  Firstly,
suppose that $B$ is a local ring and
choose a minimal reduction
$\q$ of
$\m$. Then $\ell _A(B/\q B)=2$,
because
$$\ell _A(B/\q B)=\e_{\q} ^0(B)
=\e_{\q} ^0(A)=\e(A).$$ Consequently $[
B/\n : A/\m ]{\cdot}\ell _B (B/\q
B)=2$. If
$A/\m
\ne B/\n$, then $\ell _B (B/\q
B)=1$ so that $\n =\q B$. Hence
$B$ is a regular local ring with $\n =\m$ and the assertion
follows from Proposition (2.7).
If $A/\m =B/\n$, then $\ell _B(B/\q B)=2$. Hence
$\ell_B(\n/\q B) = 1$ and so $B$ is
a Gorenstein ring, because the ring
$B/\q B$ is Artinian and Gorenstein.
The assertion now follows from Theorem
(2.3), since $\ell _A(B/A) = 1$.

Assume that $B$ is not a local
ring. Then the proof of [G2,
Proposition 4.2] still works in our
context to show that $A$ is a reduced
ring with
$\sharp\Ass A =2$, say $\Ass A=\{ \p
_1,
\p _2 \}$, such that $\p _1 +\p _2
=\m$, and $A/\p _i$ is a regular
local ring with $d = \dim A/\p _i$ for
$i=1, 2$. Hence the equality $I^2=QI$ 
follows from Theorem (3.1). This
completes the proof of Theorem
(1.1).  
\qed
\enddemo

\head 4. Examples. \endhead

In this section we shall construct 
examples, showing that the
equality $I^2=QI$ fails in general to hold, even though $A$
is a Buchsbaum local ring with sufficiently large depth and
multiplicity.

Let $k$ be a field and let $1 \leq d<m$ be integers. Let
$S=k[X_1,\cdots , X_m,V,A_1,\cdots ,A_d]$ be the polynomial
ring with $m+d+1$ variables over
$k$ and put 
$$ \frak{a} =( X_i \mid 1 \leq i \leq m-1 )^2 +( X_m ^2 )+(
X_iV
\mid 1 \leq i
\leq m ) +( V^2 -\sum _{i=1}^d A_i X_i ).$$
We regard S as a $\Bbb Z$-graded ring such that $S_0=k$ and
$X_i, V, A_j \in S_1$ for all
$1
\leq i \leq m$ and $1 \leq j \leq d$.
We put
$$R=S/\frak{a}, M=R_+, B=R_M,    \
\text{and}  \
\ \m =MB.$$ Then
$\dim R=d$, because $\sqrt{\frak{a}}
=( X_i \mid 1\leq i \leq m )+( V )$.
Let $x_i, v$, and $a_j$ denote
respectively the reduction of $X_i,
V$, and $A_j$ mod $\frak{a}$. We put $Q=(
a_1,a_2,\cdots ,a_d)$ and $\p =(
x_1,x_2,\cdots ,x_m ) +( v )$. Hence
$M=Q+\p$ and $M^3=Q\p ^2$, since $\p ^3=(0)$. Consequently,
$\q = QB$ is a minimal reduction of the maximal ideal $\m$
in $B$.

Let us begin with the following.

\proclaim{Lemma (4.1)}
$\ell _B(B/\q )=2m+1$, $\roman{e}_{\m} ^0(B)=2m$, and $B$ is
not a Cohen-Macaulay ring.
\endproclaim

\demo{Proof}
Since $R/Q \cong
k[X_1,X_2,\cdots,X_m, V]/\frak{b}$
where
$$\frak{b} =( X_i \mid 1 \leq i 
\leq m-1 )^2 +(X_m^2 , V^2 )+
(X_iV \mid 1 \leq i \leq m),$$
we have
$\dim_kR/Q=2m+1$ and $$Q : M
=Q+(x_ix_m \mid 1 \leq i \leq m-1 )+( v
).$$ Hence $\ell _B (B/\q)=2m+1$. We
put $P=( X_i \mid 1\leq i\leq m )+ (
V )$; hence $\p =P/\frak{a}$. Then
$\Min B=\{ \p \}$ and $B/\p$ is a 
regular local ring, so that we have 
$\roman{e}(B)=\ell _{B_{\p}}(B_{\p})$.
Let $\widetilde S =S[\frac{1}{A_1}]$
and $\widetilde k =
k[A_1,
\frac{1}{A_1}]$.  Then $$\widetilde S
=\widetilde k [X_1',\cdots ,X_m', V',
A_2',\cdots ,A_d'],$$ where
$X_i'=\frac{X_i}{A_1}$, $V'=\frac{V}{A_1}$, and $A_j'=\frac{A_j}{A_1}$ (
$1\leq i\leq m$, $1\leq j \leq d$ ).
The elements $\{ X_i'\} _{1\leq i\leq m}$, $V'$, and $\{ A_j'\} _{2\leq
j\leq d}$ are algebraically 
independent over $\widetilde k$. We
have $$\frak{a} \widetilde S =(X_i'
\mid 1\leq i\leq
m-1)^2+({X_m'}^2)+(X_i'V' \mid 1\leq
i \leq m)+({V'}^2-\sum _{i=1}^d
A_i'X_i')$$ and $P\widetilde S =(X_i'
\mid 1\leq i\leq m)+(V').$ Hence
$X_1'-({V'}^2-\sum _{i=2}^d
A_i'X_i') \in \frak{a} \widetilde
S$. Let $T=\widetilde{k}[X'_{2} ,
\cdots , X'_{m}, V', A_{2}', \cdots ,
A_{d}']$ and we identify  $T =
\widetilde S/(X_{1}'-
({V'}^2-\sum_{i=2}^{d} A_{i}'
X_{i}'))$. Then 
$$ \align 
\frak{a}  T 
 & =  ( {V'}^{2} - \sum_{i=2}^{d} 
 A_{i}'X_{i}', \{X_{i}'\}_{2
\leq i
\leq m-1})^{2}  + ( {X_{m}'}^{2}) +
\left(({V'}^2
 -
 \sum_{i=2}^{d}A_{i}'X_{i}'
){\cdot} V'\right) 
 \\
  &  \ \ \  +   ( X_{i}' V' \mid 2
\le i
\le m )
\\
  & =  ( {V'}^2, \{X_{i}'\}_{2 \leq
 i
\leq m-1})^{2}
  +( {X_{m}'}^2 ) + 
     ( {V'}^3 ) + ( X_{i}'V' \mid 2
\le i \le m ) 
     \ \ \  ( \text{since} \ d < m
) \\
  &=  ( X_{i}' \mid 2 \le i \le m-1 )^2
+ ( {X_{m}'}^2 ) + ( {V'}^3 ) + 
       ( X_{i}' V' \mid 2 \le i \le m
),
\endalign
$$ and $PT  
 = ( X_{i}' \mid 2 \le i \le m ) + (
V' ).$ Therefore 
$ \ell_{B_{\frak p}}(B_{\frak p}) = \ell_{S_p } (S_p/\frak a
S_p) =
\ell _U (U)$, where $$
\align
U&=k(A_{1}'
, A_{2}' , \cdots , A_{d}')
[X_{2}',
\cdots , X_{m}' , {V'}
]/\frak{b} \ \ \ \text{and}\\
\frak{b} &= (X_{i}' \mid 2 \le i \le
m-1 )^2  + ( {X_{m}'}^2 ) + ( {V'}^3 ) +
(X_{i}' {V'} \mid 2 \le i \le m ).
\endalign
$$ 
Consequently, $ \e(B) = \ell_{U} (U) =
2m < \ell_B (B / \frak q ) = 2m+1$, whence $B$ is not a
Cohen-Macaulay ring.
\qed
\enddemo

Let $\underline{a}^2 $ denote the 
sequence $ a_{1}^{2} , a_{2}^{2} ,
\cdots , a_{d}^{2} $ and let $\e^{0}_{(\underline{a}^2) B} (B) $ denote
the multiplicity of $B$ with respect to the parameter ideal
$(\underline{a}^2) B =(a_{1}^{2} , a_{2}^{2} ,
\cdots , a_{d}^{2})B$. We then have the following. 

\proclaim{Proposition (4.2)}
$ \ell_B (B / (\underline{a}^2) B ) - \e^{0}_{(\underline{a}^2) B} (B)
=1.$
\endproclaim

\demo{Proof}
Since $B$ is not a Cohen-Macaulay ring, $\ell_{B} (B / (\underline{a}^2)
B ) - 
\e ^{0}_{(\underline{a}^2) B} (B) > 0 $. It 
suffices to show $\ell_{B} (B /
(\underline{a}^2) B ) - 
\e ^{0}_{(\underline{a}^2) B} (B) \le 1 $. 
Let $$\frak{c} = (A_{i}^2
\mid 1 \le i \le d ) + (X_{i}^2 \mid 1
\le i \le m) + (V^2 - \sum_{i=1}^d A_i
X_i)$$ and put $C = S/\frak{c}
$ and 
$D= S/\left( \frak a + (A_i^2 \mid 1
\le i \le d ) \right) $.  Then $D$ is
a homomorphic image of $C$. The ring
$C$ is a complete intersection with
$\dim_kC = 2^{d+m+1}$. Let
$x_i , v,
$ and $a_j$ denote, for the moment,
the reduction of $X_i, V, $ and $A_j$
mod
$\frak a + (A_i^2 \mid 1 \le i \le d
)$. Let $
\Lambda = \{1, 2, \cdots d \} $ and
$\Gamma = \{1, 2, \cdots , m \} $. For given subsets $I \subseteq
\Lambda
$ and $J \subseteq \Gamma $ we put 
$$a_I = \underset{i \in I }\to{\Pi}
a_i
 \ \ \text{and} \ \  x_J =
\underset{j
\in J  }\to{\Pi}x_j.$$
Then the elements $\{ a_I x_J \}_{I
\subset
\Lambda, J \subset \Gamma }$ and
$\{a_I x_J v \}_{I \subset \Lambda, J
\subset \Gamma }$ span the $k$-space
$D$, because their preimages in $C$
form a $k$-bases of $C$. Notice that
$a_I x_J v = 0$ if $J \neq \emptyset $
and that $J \subseteq \{i, m \} $ for
some
$1 \le i \le m-1 $ if $x_J \neq 0 $. Hence the $k$-space $D$ is actually
spanned by the following $2^{d} (2m+1)$
elements 
$$  
a_I, \ x_ia_I, \ x_ma_I,\ 
x_ix_ma_I,  \  \text{and} \  \ a_Iv
\
\ 
\text{with}
\
\  I
\subseteq
\Lambda , 1 \le i \le m-1
. \tag 4.3
$$

Let $1 \le i \le d $ and 
$K \subseteq \Lambda $. Assume
that $i
\notin K $ but $\{ 1, \cdots , i-1 \} \subseteq K $. 
Then since
$(\sum_{i=1}^da_ix_i)(x_ma_K
)=v^2x_ma_K = 0$, we have
$$\sum_{1 \le j < i } (a_jx_j
)(x_ma_K )+ (a_ix_i )(x_m a_K)
+ \sum_{i<j \le d } (a_jx_j )(x_ma_K
)=0. \tag 4.4$$ Notice that
$\sum_{1
\le j < i } (a_jx_j )(x_ma_K ) = 0,$
since $(a_\ell x_\ell)(x_ma_K ) =
(a_\ell a_K)(x_\ell x_m )=0$ for
all $\ell
\in K$. If $i < j \leq d$ and
$j
\not\in K$, then $(a_jx_j)(x_ma_K) =
(x_jx_m)a_{K\cup \{ j
\}}$. Consequently by (4.4) we have
the following expression 
$$x_ix_ma_{K \cup \{ i \}} = (a_i
x_i)(x_ma_K ) = - \sum_{i < j \le
d, j \notin K } (x_jx_m )a_{K \cup
\{ j \} } 
$$
of $x_ix_ma_{K \cup \{ i \}}$.
Hence, letting
$K = I \setminus \{i\}$, it follows
from this expression that
for all
$1
\leq i \leq d$, the
set $\left\{x_ix_ma_I\right\}_{
\{ 1,
\cdots , i \} \subseteq I
\subseteq \Lambda}$ is contained in
the $k$-subspace of $D$ spanned by
$\{ x_j x_m a_J \}_{i < j \le d, J 
\subseteq \Lambda }$. Therefore, in
order to span the whole $k$-space $D$,
for each
$1
\le i
\le d
$ the elements
$\{x_ix_ma_I \mid \{1, 2, \cdots, i\}
\subseteq I
\subseteq \Lambda\}$ can be deleted
from the system of generators given by
(4.3), so that we have 
$$
\align
\ell_B(B/(\underline{a}^2)B) =\dim_k
D &\le 2^d (2m+1)-\sum_{i=1}^d
2^{d-i} \\ &=2^d (2m+1)-(2^d -1) \\
&=2^{d+1} m +1
\\
&=2^d\roman{e}^{0}_{(\underline{a})
B}(B) + 1 \\
&=\roman{e}^{0}_{(\underline{a}^2)
B}(B) + 1
\endalign
$$
as is claimed. 
\qed
\enddemo

By Proposition (4.2) we get $$\ell_B
(B/(\underline{a}^2)B) -
\e^{0}_{(\underline{a}^2) B}(B) =
\ell_B (B/\q) -
\e^{0}_{\q} (B) =1.$$
Hence the  local cohomology modules
$\H^i_{\frak m } (B)  ~ (i \ne d ) $
are finitely generated $B$-modules
and
$\sum_{i=0}^{d-1} {d-1\choose i}
h^i(B)=1$ (cf. \cite{SV, Appendix,
Theorem and Definition 17}).
Accordingly, either
$\depth B =0$, or $\depth B = d-1 $.
We have that $h^0 (B) =1$ and $h^i
(B) =0\ \ (1 \le i \le d-1)
$ if $\depth B =0$, and that
$h^{d-1}(B)=1$ if $\depth B = d-1$.
In any case, $\H^i_{\frak m}(B)=(0)$
for all $ i \ne t, d $, and
$\H^t_{\frak m}(B) \cong B/\frak m$,
where $t=\depth B$.  Thus $B$ is a
Buchsbaum ring (cf. \cite{SV,
Chap.~I, Proposition 2.6}). We
actually have the following. 

\proclaim{Theorem (4.5)}
$\H^{d-1}_{M}(R)\cong (R/M)(d-3). $
\endproclaim

\demo{Proof} (1) ($d=1$). Use the
fact that
$0
\ne x_1 x_m \in (0): M $.

(2) ($d=2$). Assume that $\depth B=0
$. Then applying functors
$\H^i_{M}(*)$ to the exact
sequence 
$$0 \to \H^0_{M}(R) \to R(-1) 
\overset{a_2}\to{\to} R \to R/a_2 R
\to 0, \tag 4.6$$ we get a natural
isomorphism 
$\H^0_{M}(R) \cong
\H^0_{M }(R/a_2 R ).$
We apply the result of the case where 
$d=1$ to the ring $R/a_2 R$
and choose 
$0 \ne \varphi \in R_2$ such that 
$M \varphi =(0) $ and $\varphi \equiv
x_1x_m$ mod $a_2R.$ Let 
$\varphi = x_1 x_m + a_2 \psi$ with
$\psi \in R_1$.  Then $x_1x_m a_2 +
a_2^2
\psi =0$,  that is 
$X_1X_m A_2 + A_2^2\xi \  \in \frak a
$ for some  $\xi \in S_1$, which is
impossible. Hence 
$\depth B =1 $ and by (4.6) we get 
an isomorphism $\H^0_{M}(R/a_2R)
\cong \H^1_{M }(R)(-1)$. Thus
$\H^1_{M }(R) \cong (R/M)(-1)$,
because $\H^0_{M}(R/a_2R)
\cong (R/M)(-2)$.

(3) ($d \geq 3$). We may assume that 
our assertion holds true for
$d-1$. Then $\depth B = d-1$, because
$$h^0(B/a_dB) = h^0(B)+ h^1(B)$$ and
$h^0(B/a_d B) = 0$ by the 
hypothesis on $d$. Hence $a_d$ is a
non-zerodivisor of $R$, so that we
have
$$\H_M^{d-2}(R/a_d R) \cong 
\H_M^{d-1}(R)(-1).$$ Thus
$\H_M^{d-1}(R)
\cong (R/M)(d-3)$, because
$\H_M^{d-2}(R/a_dR)\cong (R/M)(d -
4)$.
\qed
\enddemo

Let $J=Q:M$ and $I= \q:\m$~$(= JB)$.
We then have the following.

\proclaim{Proposition (4.7)}
$I^2 \ne \q I$ but $I^3 = \q I^2$.
\endproclaim

\demo{Proof}
We have $J = Q + (x_i x_m
\mid 1 \leq i \leq m-1) + (v)$ (cf.
Proof of Lemma (4.1)). Assume that 
$J^2 = QJ$. Then $v^2 \in QJ$ whence
$$V^2 \in (A_i \mid 1 \leq i \leq d) 
{\cdot}\left[(A_i \mid 1 \leq i \leq
d) + (X_i X_m \mid 1 \leq i \leq m-1)
+ (V)\right] + \frak{a}. \tag 4.8$$
We now substitute $X_i = 0$
and $A_j = 0$ for all $2 \leq i \leq
m$ and $2
\leq j \leq d$. Then by (4.8) we get
$$V^2 \in (A_1^2, A_1V, X_1^2, X_1V,
V^2-A_1X_1)$$
in the polynomial
ring $k[X_1,V,A_1]$, which is
impossible. Hence $v^2 \not\in QJ$ so
that we have $J^2 \ne QJ$. We get $J^3
= QJ^2$, because $J^2 = QJ + (v^2)$
and
$v^3 = 0$. 
\qed
\enddemo

Therefore,
for given integers $1 \leq d < m$,
there exists a Buchsbaum local ring 
$A$ with $\dim A = d$, 
$\depth A = d-1$, and $\roman{e}(A)=
2m$, such that
$A$ contains a parameter ideal $Q$
which is a minimal reduction of $\m$
and $I^2 \ne QI$ but $I^3 = QI^2$,
where $I = Q:\m$.  Thus the equality
$I^2 = QI $ fails in general to hold,
even though $A$ is a Buchsbaum local
ring with sufficiently large depth
and multiplicity.

%%%%%%%%%%%%%%%%%%%%%%%%%  
%%%%%%ref   
%%%%%%%%%%%%%%%%%%%%%%%%%%%%%%%%%%%


\Refs

\widestnumber\key{CHV}


\ref
\key CHV
\by A. Corso, C. Huneke, and W. V. Vasconcelos
\paper On the integral closure of ideals
\yr 1998
\pages 331-347
\vol  95
\jour manuscripta math.
\endref

\ref
\key CP
\by A. Corso and C. Polini
\paper Links of prime ideals and their Rees algebras
\yr 1995
\pages 224-238
\vol  178
\jour J. Alg.
\endref





\ref
\key CPV
\by A. Corso,  C. Polini, and W. V. Vasconcelos
\paper Links of prime ideals
\yr 1994
\pages 431-436
\vol  115
\jour Math. Proc. Camb. Phil. Soc.
\endref





\ref
\key    G1
\by     S.~Goto
\paper  On the Cohen-Macaulayfication 
of certain Buchsbaum rings
\jour   Nagoya Math. J.
\vol    80
\yr     1980
\pages  107-116
\endref

\ref
\key    G2
\by     S.~Goto
\paper   Buchsbaum rings 
with multiplicity 2
\jour   J. Alg.
\vol    74
\yr     1982
\pages  494-508
\endref


\ref
\key G3
\by S.~Goto
\paper Integral closedness of complete-intersection ideals
\yr 1987
\pages 151-160
\vol  108
\jour J. Alg.
\endref

\ref
\key GH
\by S. Goto and F. Hayasaka
\paper Finite homological dimension and
primes associated to integrally closed
ideals II 
%\yr
%\pages 43-54
%\vol  135
\jour J. Math. Kyoto Univ. 
\toappear
\endref

\ref
\key GN
\by S.~Goto and K.~Nishida
\paper Hilbert coefficients and
Buchsbaumness of associated graded rings
\yr 2003
\pages 61-74
\vol  181
\jour J. Pure and Applied Alg.
\endref



\ref
\key GSa
\by S.~Goto and H.~Sakurai
\paper The equality $I^2=QI$ in Buchsbaum
 rings
%\yr
%\pages
%\vol 
\jour Rendiconti del Seminario Matematico dell'Universit di Padova (to appear)
\endref

\ref
\key GSu
\by S.~Goto and N.~Suzuki
\paper Index of reducibility of
parameter ideals in a local ring
\yr 1984
\pages 53-88
\vol  87
\jour J. Alg.
\endref

\ref
\key    SV
\by     J. St\"{u}ckrad and W. Vogel
\book   Buchsbaum rings and applications
\publ Springer-Verlag
\publaddr  Berlin, New York, Tokyo
\yr 1986
\endref

\ref
\key    Y1
\by     K. Yamagishi
\paper  The associated graded modules of 
Buchsbaum modules with respect to 
$\frak{m}$-primary ideals in the
equi-$\Bbb{I}$-invariant case
\jour   J. Alg.
\vol    225
\yr     2000
\pages  1-27
\endref

\ref
\key    Y2
\by     K. Yamagishi
\paper  Buchsbaumness in Rees modules associated to ideals of
minimal multiplicity in the equi-$\Bbb
I$-invariant case
\jour   J. Alg.
\vol    251
\yr     2002
\pages  213-255
\endref




\endRefs
\enddocument